\documentclass[12pt,leqno]{article}
\usepackage{amssymb,amsmath,amsthm}
\usepackage[all]{xy}
\usepackage[alphabetic,lite]{amsrefs}
\numberwithin{equation}{section}
\usepackage{mathrsfs}
\usepackage[top=2cm,bottom=2cm,left=3cm,right=3cm,marginparsep=0.3cm,marginparwidth=3cm,includefoot]{geometry}
\usepackage{graphicx}
\usepackage{color}
\newcommand{\scbul}{{\,\raise.4ex\hbox{$\scriptscriptstyle\bullet$}\,}}

\newcommand{\R}{\mathbb{R}}

\date{}
\begin{document}

\title{Mikio Sato, a visionary of mathematics
}

\begin{figure}\centerline{
\begin{tabular}{cc}
\includegraphics[scale=.5437]{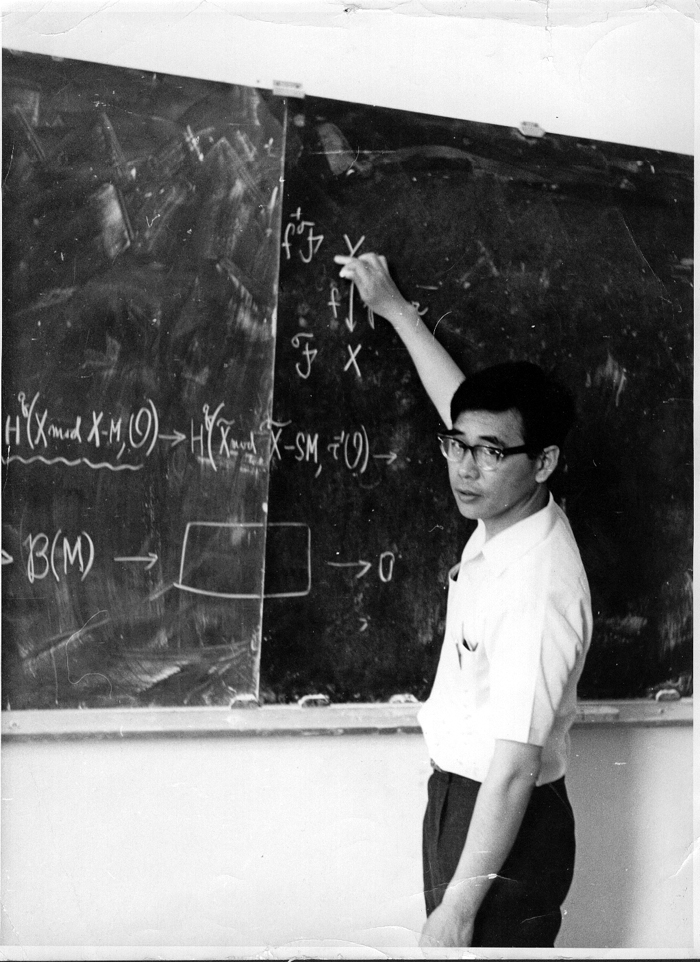}
\end{tabular}
}\caption{Mikio Sato, Nice (France) 1970}
\end{figure}

\author{Pierre Schapira\footnote{The part of this paper concerned with number theory   had benefitted from the scientific comments of Jean-Beno{\^i}t Bost and Antoine Chambert-Loir, and also Pierre Colmez for the last version.  I warmly thank all of them.}
}

\begin{abstract}
This paper, to appear in the  ``Notices of the AMS'' 2024, is a modified version of a text already appeared 
in this journal, Feb. 2007 after a first publication 
in French, in  ``La Gazette des
Math{\'e}maticiens'' {\bf 97} (2003) on the occasion of Sato's 
reception of the 2002/2003 Wolf prize. 
\end{abstract}
\maketitle

\subsection*{}
Mikio Sato passed away on  January 9, 2023 and it was very sad news for all of us  who had the chance to meet him and share his vision of mathematics. A  vision radically new and revolutionary, too revolutionary to be immediately understood by the community of mathematicians, in particular  by the analysts, despite the fact that initially, Sato's aim was to develop ``algebraic analysis'', that is, to treat problems of analysis with the tools of algebraic geometry. 
 And Sato never did any effort nor spend a great deal of time or energy 
  to make his ideas propagate. 
  
Sato did not write a lot, did not communicate easily and attended 
very few meetings. But he
invented a new way of doing analysis,  ``Algebraic Analysis'', and he opened a new horizon in mathematics, the microlocal approach. Thanks to Sato, we understand that many phenomena which appear on a manifold are in fact the projection on the manifold of things living in the cotangent bundle to the manifold. Being ``local''
becomes, in some sense, global with respect to this projection. 

Mikio Sato  also  created a school, ``the Kyoto school'', among whom Masaki Kashiwara, Takahiro Kawai, 
Tetsuji Miwa and Michio Jimbo should be mentioned. 

Born in 1928\footnote{See~\cite{An07} for more details about Sato's life.},
Sato became known in mathematics only in 
1959--60 with his theory of hyperfunctions.
Indeed, his studies had been seriously disrupted by the war,
particularly by the bombing of Tokyo. 
After his family home was burned down, he had to work as a coal delivery
  man and later as a school teacher.  At age 29 he became an assistant 
  at  the University of Tokyo. He studied mathematics and
physics, on his own.

To understand the originality of Sato's theory of hyperfunctions, one
has to place it in the mathematical landscape of the time.
Mathematical analysis from the 1950s to the 1970s was under the
domination of functional analysis, marked by the success of the theory
of distributions. People were essentially looking for existence
theorems for linear partial differential equations (LPDE) and most of the
proofs  were reduced to finding ``the right functional space'',
to obtain some  a priori estimate and apply the  Hahn--Banach theorem.

It was in this environment that Sato defined hyperfunctions~\cite{Sa60} in 1959--1960 
as boundary values of holomorphic functions, a discovery 
which allowed him to obtain a position at the University of Tokyo 
thanks to the clever patronage of Professor Iyanaga, an exceptionally
open-minded person and a great friend of French culture. 
Next, Sato spent two years in the USA, in Princeton, where
he unsuccessfully tried to convince  Andr{\'e} Weil of the relevance 
of his cohomological approach to analysis.

Sato's method was radically new, in no way using the notion of limit. 
His hyperfunctions are not limits of functions in any sense of the
word, and the space of hyperfunctions has no natural topology other than the
trivial one. For his construction, Sato invented local cohomology 
in parallel with Alexander Grothendieck.
This was truly a revolutionary vision of analysis. 
And, besides its evident originality, Sato's approach 
had deep implications since it naturally led to microlocal
analysis. 

The theory of LPDE with 
variable coefficients was at its early beginnings in the years 1965--1970
and under the shock of  Hans Lewy's example which showed that 
a very simple first order linear equation                                               
$(-\sqrt{-1}\partial_1 +\partial_2-2(x_1+\sqrt{-1}x_2)\partial_3)u=v$
had no solution, even a local solution, in the space of distributions\footnote{The slightly simpler equation
$(\partial_1 +\sqrt{-1}x_1\partial_2)u=v$
does not have any solution in the space of germs at the origin
of distributions in  $\R^2$ either, nor even in the space of germs of
hyperfunctions.}.
The fact that an equation had no solution was quite disturbing at that
time. People thought it was a defect of the theory, that the spaces
one had considered were too small to admit solutions. 
Of course, often just the opposite is true and one finds that the
occurrence of a cohomological obstruction heralds interesting
phenomena:
the lack of a solution is the 
demonstration of some deep and  hidden geometrical phenomena. In the case
of the  Hans Lewy equation, the hidden geometry is
``microlocal'' and this equation is microlocally equivalent to an
induced  Cauchy--Riemann equation on a real hypersurface of the complex space.

\begin{figure}\centerline{
\begin{tabular}{cc}
\includegraphics[scale=.7437]{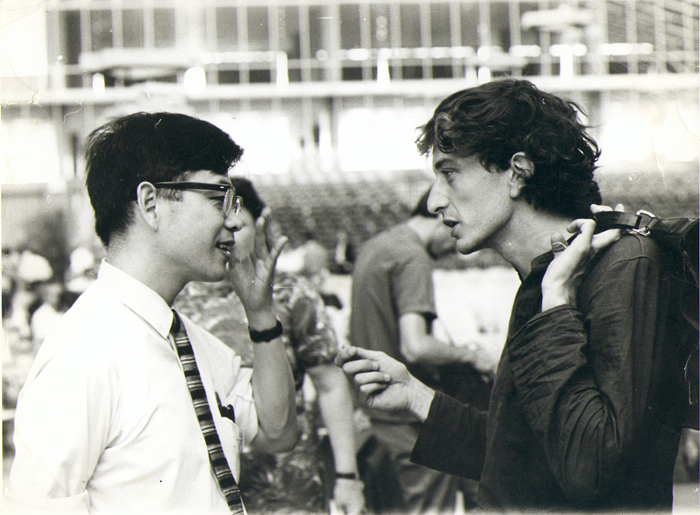}
\end{tabular}
}\caption{Mikio Sato and the author, ICM, Nice 1970}
\end{figure}

In mathematics, as in physics, in order to treat phenomena in a given
(affine) space, one is naturally led to compute in the dual
space. One way to pass from a vector space to his dual, the most commonly used in analysis, is the Fourier
transform. This transform, being not of local nature, is not easily adapted to calculus on manifolds.
By contrast, Sato's
method is perfectly suited to manifolds. If $M$ is a real analytic manifold, $X$ a complexification of $M$,
what plays the role of the dual space is now the conormal bundle $T^*_MX$ to $M$ in the cotangent bundle $T^*X$, something local on $M$. 
(Note that $T^*_MX$ is isomorphic to $\sqrt{-1}T^*M$.)
In order to  pass from $M$ to $T^*_MX$, Sato 
constructed a key tool of sheaf theory, the microlocalization functor $\mu_M(\scbul)$, 
the ``Fourier--Sato'' transform of the specialization functor 
$\nu_M(\scbul)$.
This is how Sato defines in~\cite{Sa70}  the analytic wave front set of
hyperfunctions (in particular, of distributions), a closed conic 
subset of the cotangent bundle, and 
he shows that if a hyperfunction $u$ 
is a solution of the equation $Pu=0$, then its  wave front set
is contained in the intersection with $T^*_MX$ of the characteristic variety of the
operator $P$. It is then clear (but it was not so clear at this time) that if you want to understand what happens on a real manifold, you better look at what happens on a complex neighborhood of the manifold. 

Of course, at this time other mathematicians (especially
Lars H{\"o}rmander) and physicists ({\em e.g.,} Daniel Iagolnitzer) had the intuition
that the cotangent bundle was the natural space for analysis, and in
fact this intuition arose much earlier,  in particular in the work of Jacques Hadamard,
Fritz John and Jean Leray.
Indeed, pseudo-differential operators did
exist before the wave front set.  

In 1973, Sato and his two students, M.~Kashiwara and T.~Kawai, 
published a treatise~\cite{SKK73} on the microlocal analysis of LPDE. 
Certainly this work had a considerable impact, although most 
analysts didn't understand a single word. H{\"o}rmander and his school 
then adapted the classical Fourier transform 
to these new ideas, leading to the now popular theory of
Fourier-integral operators (see for example~\cites{Ho83}). 

The microlocalization functor is the  starting point of microlocal analysis but it is   
 also at the origin of the microlocal theory of sheaves, due to  Kashiwara and the author~\cite{KS90}. This theory associates to a sheaf on a real manifold $M$ its micro-support, a closed conic subset of the cotangent bundle $T^*M$ and  allows one to treat sheaves ``microlocally'' in $T^*M$. The theory of  systems of LPDE becomes essentially sheaf theory, the only analytic ingredient being the Cauchy-Kowalevsky theorem. One of the deepest result of~\cite{SKK73} was the involutivity theorem which asserts that the characteristic variety of a microdifferential system  (in particular, of a ${\cal D}$-module, see below) is co-isotropic. Similarly the micro-support of a sheaf is proved to be co-isotropic, which makes a link  between sheaf theory and symplectic topology, a link at the origin of numerous important results.
 
 Also note that by a fair return of things, microlocal analysis, through microlocal sheaf theory, appeared quite recently in algebraic geometry, under the impulse of Sacha Beilinson~\cite{Be16}.

Already in the 1960s, Sato had the intuition of 
${\cal D}$-module theory, of holonomic systems and of the $b$-function
(the so-called Bernstein--Sato $b$-function). He gave a series of
talks on these topics at Tokyo University but had to stop 
for lack of combatants. 
His ideas were reconsidered and systematically developed 
by Masaki Kashiwara in his 1969 thesis\footnote{See \cite{Sc18} for an overview of Kashiwara's work, a part of  which was deeply influenced by Sato's ideas.}.
As its name indicates,
a ${\cal D}$-module is a  module over the sheaf of rings
${\cal D}$ of differential operators, and a module over a ring
essentially means  ``a system of linear equations''\footnote{According to Mikio Sato (personal communication), at the origin of this idea is the mathematician and philosopher of the 17th century, E.~W von Tschirnhaus.} 
with coefficients in this ring. 
The task is now to treat (general) systems 
of LPDE. 
This theory, which also simultaneously appeared in the more
algebraic framework developed by Joseph Bernstein, a student of Israel  Gelfand, 
quickly had considerable success  in several branches of
mathematics. In 1970--1980, Kashiwara obtained 
almost all the fundamental results of the theory, in particular those 
concerned with holonomic modules, such as his constructibility theorem, 
his index theorem for holomorphic solutions of holonomic modules, 
the proof of the rationality of the zeroes 
of the  $b$-function and his proof of the (regular) Riemann-Hilbert 
correspondence.

The mathematical landscape of 1970--1980 had thus 
considerably changed. Not only did one treat equations with variable
coefficients, but one treated systems of such equations and moreover
one worked microlocally, that is, in the cotangent bundle, the phase
space of the physicists. But there were two schools in the world:
the $C^\infty$ school, in the continuation of  classical analysis and headed by 
H{\"o}rmander who developed the calculus of Fourier integral 
operators\footnote{Many names should be quoted at this point, in particular
 those of Viktor Maslov  and Vladimir Arnold.},
and the analytic school that Sato established, which was
 almost nonexistent outside Japan and France.

France was a strategic place to receive Sato's ideas since they 
are based on, or parallel to, those of both Jean Leray and Alexander Grothendieck.
Like Leray, Sato understood that singularities have to be sought 
in the complex domain, even for the understanding of real phenomena.
Sato's algebraic analysis is based on sheaf theory, 
a theory invented by 
Leray in 1944 when he was a prisoner of war, 
clarified by Henri Cartan and made extraordinarily efficient by
 Grothendieck and his formalism of derived categories and the 
{\em six op{\'e}rations}.

Sato, motivated by physics as usual, then tackled the analysis of the
$S$-matrix in light of microlocal analysis.  With his two
new students, M.~Jimbo and T.~Miwa~\cite{SMJ78}, he
explicitly constructed the solution of
the $n$-points function of the Ising model in dimension $2$ using 
Schlesinger classical theory of isomonodromic deformations of ordinary
differential equations. This naturally led him to the study of
KdV-type non-linear 
equations. In 1981, with his wife Yasuko Sato (see~\cite{SaSa82} and also~\cite{Sa89}), he 
interpreted the solutions of the KP-hierarchies as points of an
infinite Grasmannian manifold and introduced his famous $\tau$-function.
These results would be applied to other classes of equations and 
would have a great impact in mathematical physics in the study of
integrable systems and field theory in dimension $2$. 

In parallel with his work in analysis and in mathematical physics, Sato 
obtained remarkable results in group theory and in number theory.

He introduced the theory of prehomogeneous vector spaces,
that is, of linear representations of complex reductive groups 
with a dense orbit. The important case where the complement of this
orbit is a hypersurface gives good examples of $b$-functions (see~\cite{Sa-Sh1, SK77}).

In 1962, Sato also discovered how to deduce, using a construction of auxiliary (Kuga-Sato) varieties, the 
Ramanujan conjecture on the coefficients of the modular form
$\Delta$ from Weil's conjectures  concerning the number of solutions 
of polynomial equations on finite fields.
His ideas allowed Michio Kuga and Goro Shimura to treat the case of compact
quotients of the  Poincar{\'e} half-space and one had  to wait 
another 10 years for Pierre Deligne to definitely prove 
that  Weil's conjectures imply  Ramanujan and Petersson's conjecture. 

Mikio Sato  shared the  Wolf  prize with John Tate in 2002/03.  They also share a famous
conjecture
in number theory concerning the repartition of Frobenius
angles. Let $P$ be a degree 3 polynomial with integer coefficients and
simple roots. Hasse has shown that for any prime $p$ which does not
divide the discriminant of $P$, the number of solutions of the 
congruence $y^2=P(x)\pmod p$  is like $p-a_p$,
with $\mathopen| a_p \mathclose |\leq 2\sqrt p$.
When writing $a_p= 2\sqrt p \cos\theta_p$ with $0\leq\theta_p\leq\pi$,
the Sato--Tate conjecture predicts that these angles $\theta_p$
follow the law $(2/\pi)\sin^2\theta$ (in absence of complex multiplication).
Note that Tate was led to this conjecture by the study of 
algebraic cycles and Sato by computing numerical data. 

Sato's most recent works are essentially unpublished (see however~\cite{Sa89a}) and have been 
presented in seminars attended only by a small group of people.
They treat an algebraic approach
of non-linear systems of PDE, in particular of holonomic
systems, of which theta functions are examples of solutions !

Looking back, 50 years later, we realize that Sato's approach to
mathematics is not so different from that of Grothendieck, that 
Sato did have the incredible temerity to treat analysis as algebraic
geometry and was also able to build the algebraic and geometric
tools adapted to his problems.

His influence on mathematics is, and will remain, considerable.

\providecommand{\bysame}{\stLeavevmode\hbox to3em{\hrulefill}\thinspace}

\vspace*{1cm}

{\scriptsize{
\parskip=-1ex
\noindent
Pierre Schapira\\
Sorbonne Universit{\'e}, CNRS IMJ-PRG\\
4 place Jussieu, 75252 Paris Cedex 05 France\\
e-mail: pierre.schapira@imj-prg.fr\\
http://webusers.imj-prg.fr/\textasciitilde pierre.schapira/
}}

\end{document}